# Voronoi Diagrams Generated by the Archimedes Spiral

Mark Frenkel, Irina Legchenkova, Edward Bormashenko*

Department of Chemical Engineering, Engineering Sciences Faculty, Ariel University, Ariel, Israel 40700.

* Correspondence: edward@ariel.ac.il

**Abstract:** Voronoi mosaics inspired by the seed points placed on the Archimedes Spirals are reported. Voronoi entropy was calculated for these patterns. Equidistant and non-equidistant patterns are treated. Voronoi mosaics built from cells of equal size which are of a primary importance for decorative arts are reported. The pronounced prevalence of hexagons is inherent for the patterns with an equidistant and non-equidistant distribution of points, when the distance between the seed points is of the same order of magnitude as the distance between the turns of the spiral. Penta- and heptagonal "defected" cells appeared in the Voronoi diagrams due to the finite nature of the pattern. The ordered Voronoi tessellations demonstrating the Voronoi entropy larger than 1.71, reported for the random 2D distribution of points, were revealed. The dependence of the Voronoi entropy on the total number of the seed points located on the Archimedes Spirals is reported. The aesthetic attraction of the Voronoi mosaics arising from seed points placed on the Archimedes Spirals is discussed.

**Keywords:** Archimedes Spiral; Voronoi entropy; surface patterns; aesthetic attraction; phyllotaxis.

## 1. Introduction

A spiral is a curve which emanates from a point, moving farther away as it revolves around the point. Spirals inspiring wonder and curiosity are abound in nature, mathematics, art and decoration [6-8]. A spiral like curve has been found in Mezine, Ukraine, as part of a decorative object dated to 10.000 BCE. Evergreen spiral motifs appear also on the altar found in the Temples of Malta (3000 BC) and, also, they are inherent for the Celtic megalithic culture. They are seen on Egypt Minoan pottery. The famous spirals were created by Leonardo da Vinci. Spirals are inspiring the modern artists such as Robert Smithson and Fransisco Infante-Arana. In the microscopic scale, the DNA molecules twist round in the form of two helices, whereas, in the largest possible scale, the arms of galaxies curl round in the form of logarithmic spirals [9]. The physical world exhibits a startling repetition of spiral patterns [6-8]. Biological patterns also often demonstrate spiral-like structures. In particular, geometric models of phyllotaxis were used to generate realistic images of flowers and fruits with spiral patterns [14].

In our paper we focus on patterns generated by Archimedes (or Archimedean) spiral, used for creating of the Voronoi partition. The Archimedean spiral (abbreviated for brevity AS) is a spiral with the polar equation $r = a\theta^{1/n}$, where r is the radial distance, $\theta$ the polar angle, and n is a constant which determines how tightly the spiral is ''wrapped'' [4,32]. When this constant $n = 1$, the resulting spiral is given by $r = a\theta$. In this case the AS has the property that any ray from the origin intersects successive turnings of the spiral in points with a constant separation distance (equal to 2πa if $\theta$ is measured in radians). That is why this spiral is also called the "arithmetic spiral".

The AS demonstrates numerous natural and technological exemplifications. As an example, the drawing of an AS (spirography) is commonly used in the evaluation of patients with pathologic tremors and other movement disorders [23]. Interlocked AS supplied a relief cutting method to turn rigid planar surfaces into flexible ones using meander patterns [34]. Archimedes' spiral grooves produced on silver films supplied a selective chirality to surface plasmons [25]. Artistic space-filling designs based on spiral packings were reported [4].

We exploited AS for generating Voronoi partitions, demonstrating interesting mathematical properties and aesthetic appeal. Voronoi partitions (or tessellations) enable quantification (the expression or measurement) of ordering in sets of points [2,29,31]. The idea of what is now called the Voronoi tessellation was proposed by Johannes Kepler and Rene Descartes [10,19]. Descartes used these tessellations to verify that the distribution of matter in the Universe forms vortices centered at fixed stars [10,19]. The idea was developed by Dirichlet in the context of his works on quadratic forms [11].

Let us explain the idea of the Voronoi diagram (tessellation). A tessellation or tilling the plane is a collection of plane figures that fills the plane with no overlaps and no gaps. A Voronoi diagram is a partitioning of a plane into regions (cells) based on the distance to a specified discrete set of points (called seeds, sites, nuclei or generators) [2,31]. For each seed, there is a corresponding region consisting of all points closer to that seed than to any other [2,31]. The Voronoi polyhedron of a point nucleus in space is the smallest polyhedron formed by the perpendicularly bisecting planes between a given nucleus and all the other nuclei [31]. The Voronoi tessellation divides a region into space-filling, non-overlapping convex polyhedral. Voronoi diagrams represent planar graphs [2,31]. The topological properties of Voronoi diagrams are surveyed in Ref. 31.

The Voronoi tessellation enables quantification of ordering of the 2D structure, by the calculation of the so-called Voronoi entropy defined as:

$$S_{vor} = -\sum_i P_i ln P_i, \qquad (1)$$

where $i$ is the number of polygon types, $P_i$ is the fraction of polygons possessing $n$ sides (edges) inherent for a given Voronoi diagram (also called the coordination number of the polygon) [2,12,13,29]. The Voronoi entropy becomes zero for a perfectly ordered structure (when we have just a single kind of polygons); for a typical case of fully random 2D distribution, the value of $S_{vor}$ = 1.71 was reported [20,21]. Eq. 1 has the form similar to the statistical measure of information and the entropy in statistical mechanics [28]. That is why it was called "the Voronoi entropy".

Consider some simple exemplifications of the Voronoi tessellation (see **Figure 1**). **Figure 1A** represents the regular array of points (left) that leads to a regular array of square tiles (right) with the Voronoi entropy $S_{vor}$ which equals zero (indeed $P_1 = 1; ln P_1 = 0$ in Eq. 1). **Figure 1B** represents the pattern (left) giving rise to the Voronoi tessellation built from irregular (distorted) hexagons (right); the corresponding Voronoi entropy of the tessellation, demonstrated in **Figure 1B** also equals zero (again $P_1 = 1; ln P_1 = 0$ is true for this pattern). **Figure 1C**, in turn, depicts a semi-regular set of points (left) resulting in a twin-tile tessellation (i.e., regular hexagons and smaller squares, right). The Voronoi entropy of the tessellation shown in **Figure 1C** (left) $S_{vor} = \frac{1}{3} ln \frac{1}{3} + \frac{2}{3} ln \frac{2}{3} = 0.6365$. **Figure 1D** demonstrates the pattern emerging from the 75 randomly placed points (left) and the Voronoi tessellation (right) arising from this pattern. The Voronoi entropy of this pattern $S_{vor} = 1.6959$ is close to the value $S_{vor} = 1.71$ established for the randomly distributed sets of points [20, 21]. We will demonstrate that the tessellations with $S_{vor} > 1.71$ are possible. **Figure 1E** exemplifies the regular pattern of 80 points (left), giving rise to the Voronoi tessellation (right) with an entropy larger than that inherent for randomly distributed points.

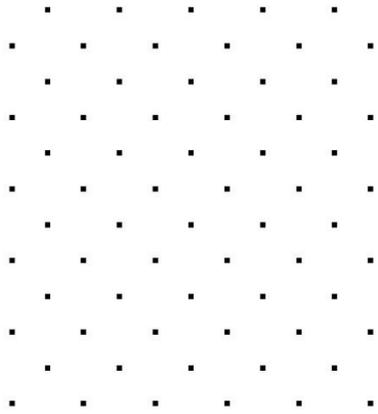
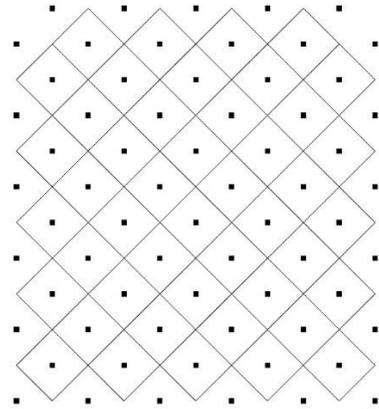

A

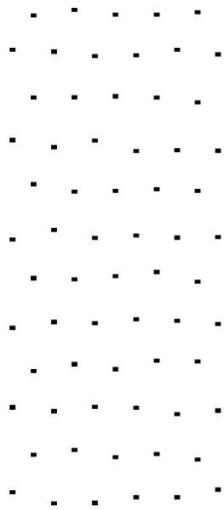
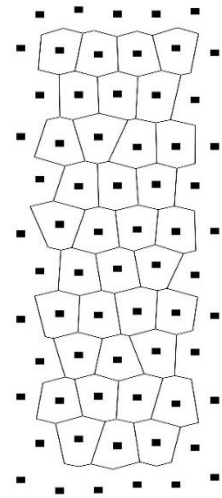

B

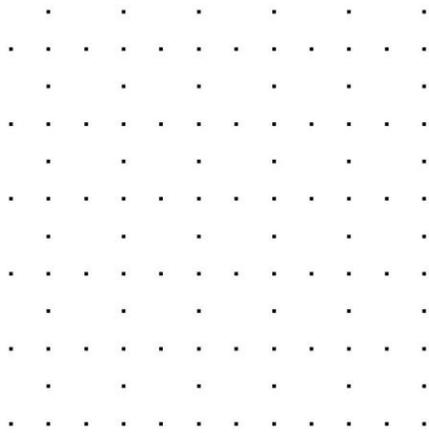
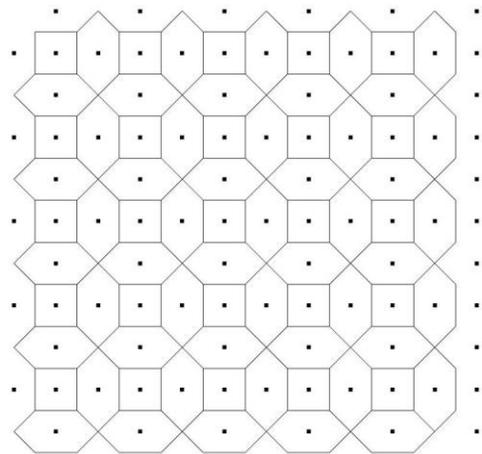

C

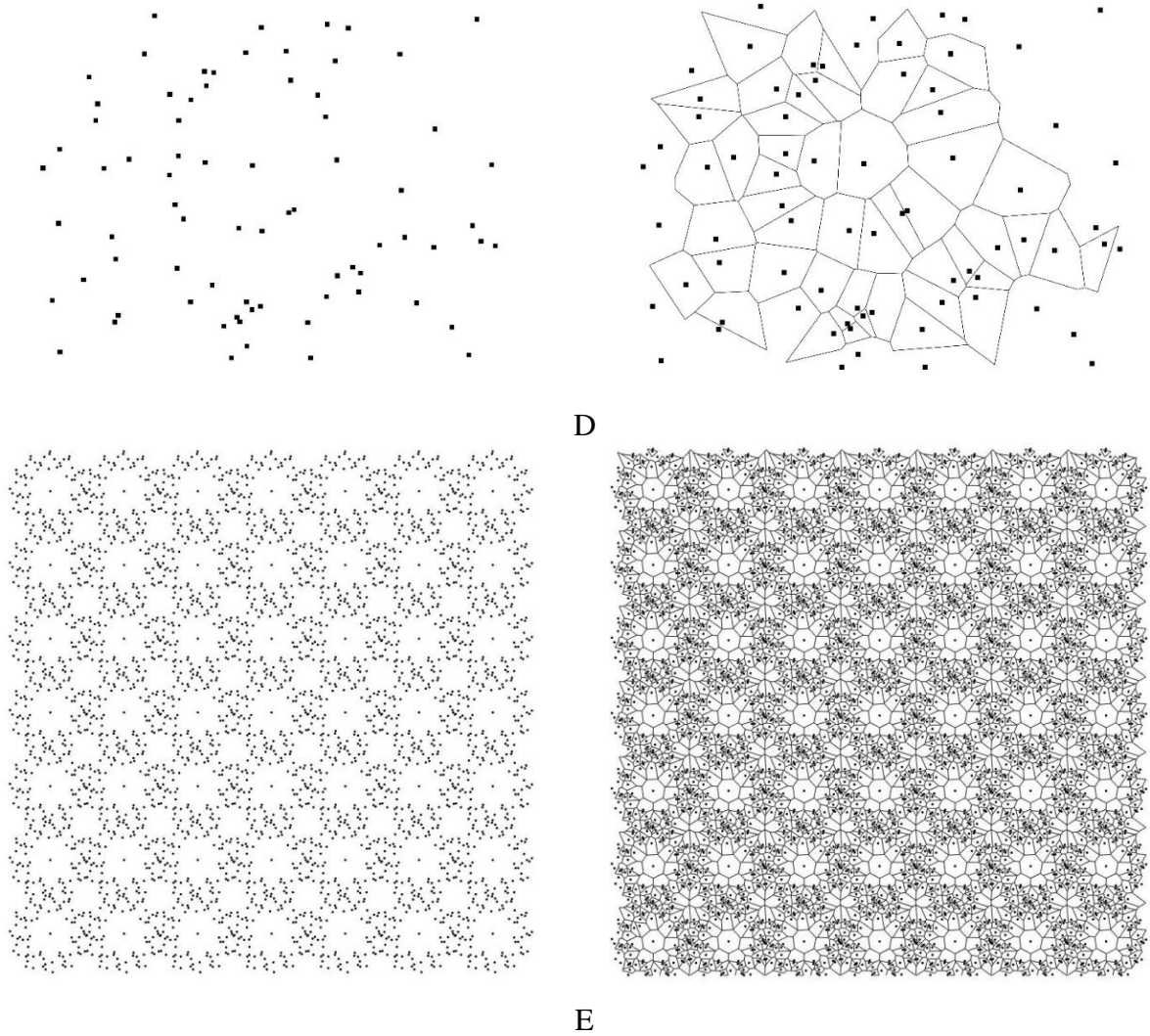

E

**Figure 1.** Exemplifications of the Voronoi tessellation. A. A regular array of points (left) that leads to a regular array of square tiles (right) with the resulting Voronoi entropy $S_{vor}$ which equals to zero. B. The pattern containing 110 points (left) giving to rise the Voronoi tessellation built of the hexagons only (right); the Voronoi entropy of the pattern is zero. C. Pattern containing 108 points (left) giving rise to the Voronoi diagram built of the regular hexagons and smaller squares (right). The Voronoi entropy of the pattern $S_{vor} = \frac{1}{3} ln \frac{1}{3} + \frac{2}{3} ln \frac{2}{3} = 0.6365$. D. The pattern emerging from 75 randomly located points (left) and the corresponding Voronoi diagram (right) are shown, the Voronoi entropy of the pattern is $S_{vor} = 1.6959$. E. Voronoi pattern arising from $7 \times 7$ *XY* translation of the pattern shown in Figure 9A is shown (left). The Voronoi tessellation (right) is built from eight types of polygons and the Voronoi entropy corresponding to the tessellation is $S_{vor} = 1.9327$.

Topological argument, arising from the Euler equation for the Voronoi diagrams is that in the limit of a large system, the average number of edges surrounding a cell is six. This leads to prevailing of hexagons in Voronoi diagrams emerging from large, random sets of points [31].

2. **Voronoi partitions generated by the AS**

The MATLAB software was used for calculation of the coordinates of the points on the AS and the subsequent generation and processing of the corresponding Voronoi patterns. To create the Voronov diagrams, we used moduli of the program developed at the Department of Physics and Astronomy at the University of California (Department of Physics and Astronomy University of California, Irvine) (https://www.physics.uci.edu/~foams/do_all.html).

AS with various parameters (points density and quantity) was generated with MATLAB software routine. Coordinates of points on AS in rectangular coordinate system were defined by the following equations:

$$x = r \cdot \cos(\varphi), \tag{2}$$

$$y = r \cdot \sin(\varphi).$$

Consider first the equidistant distribution of the points along the given AS. Parameters $r$ and $\varphi$ were varied in a way providing a given distances between the points neighboring along the spiral ($p$) and the coils of the spiral ($q$) (See **Figure A3**). The Voronoi tessellations and Voronoi entropy were established for the aforementioned points. In the case of a linear increase of distance between neighboring points (abbreviated in the text NP) on the spiral formulae (2) were transformed into Eqs. 3 (see **Appendix A,** this transformation was made in order to obtain a finite set of points and for the convenience of further calculations):

$$x_n = t_n \cdot cos(t_n),$$

$$y_n = t_n \cdot sin(t_n), \tag{3}$$

here $x_n$ and $y_n$ are the coordinates of a single point on a spiral. The variable $t$ is an array of values that alter discretely from $b$ to $d$ with a step of $c$. Parameters $b$, $c$, and $d$ allow obtaining a finite set of coordinates for developing different AS. Parameter $b$ sets a value of spiral starting (origin) point coordinates. For a sake of simplicity, we adopt $b=0$ for all of the studied patterns; this assumption corresponds to the spirals starting from the coordinate origin. Modification of the parameters $c$ and $d$ enabled generation of points located on the AS with the controlled distances between them. Aforementioned parameters $p$ and $q$ denoting linear dimensions were taken in millimeters, parameters $b$, $c$, $d$, and $t$ were dimensionless. Note, that from the "physical point" of view the dimensionless parameter $\xi = \frac{p}{q}$ appears. The ratio $\frac{p}{q}$ controls shape of the spiral and defines the distribution of points on it. Thus, the value of $\xi$ influences the properties of the Voronoi tessellation, as it will be demonstrated below.

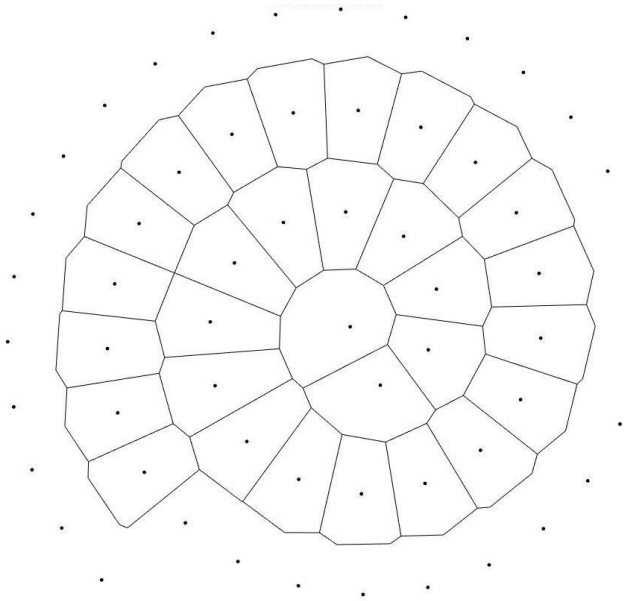

A

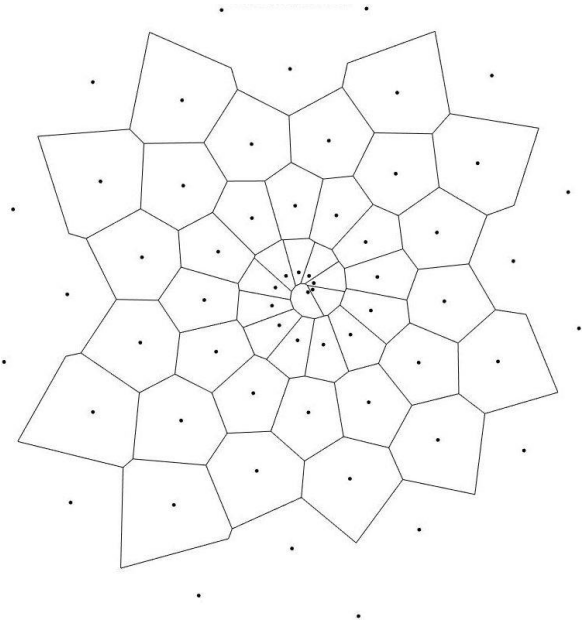

B

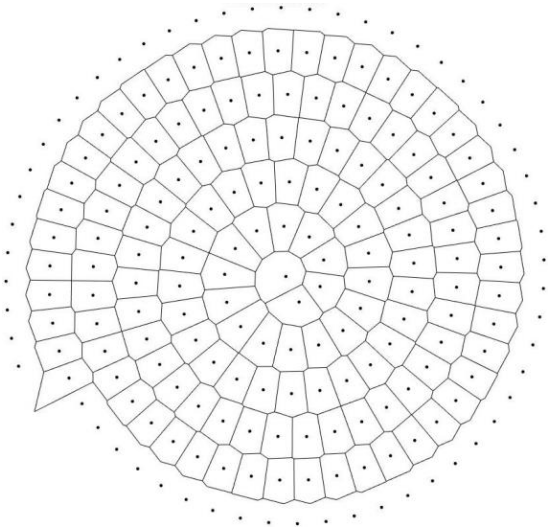

C

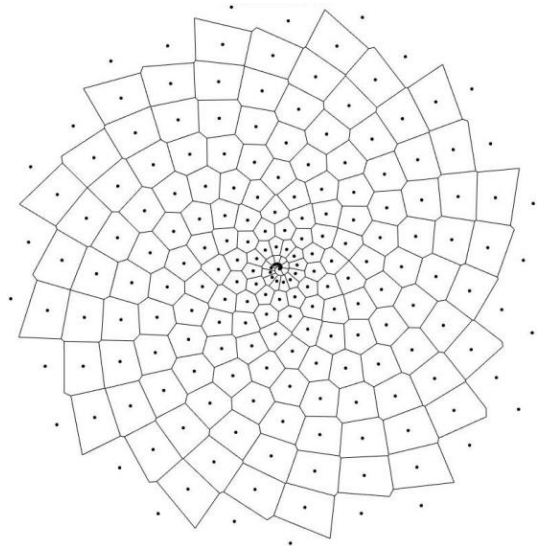

D

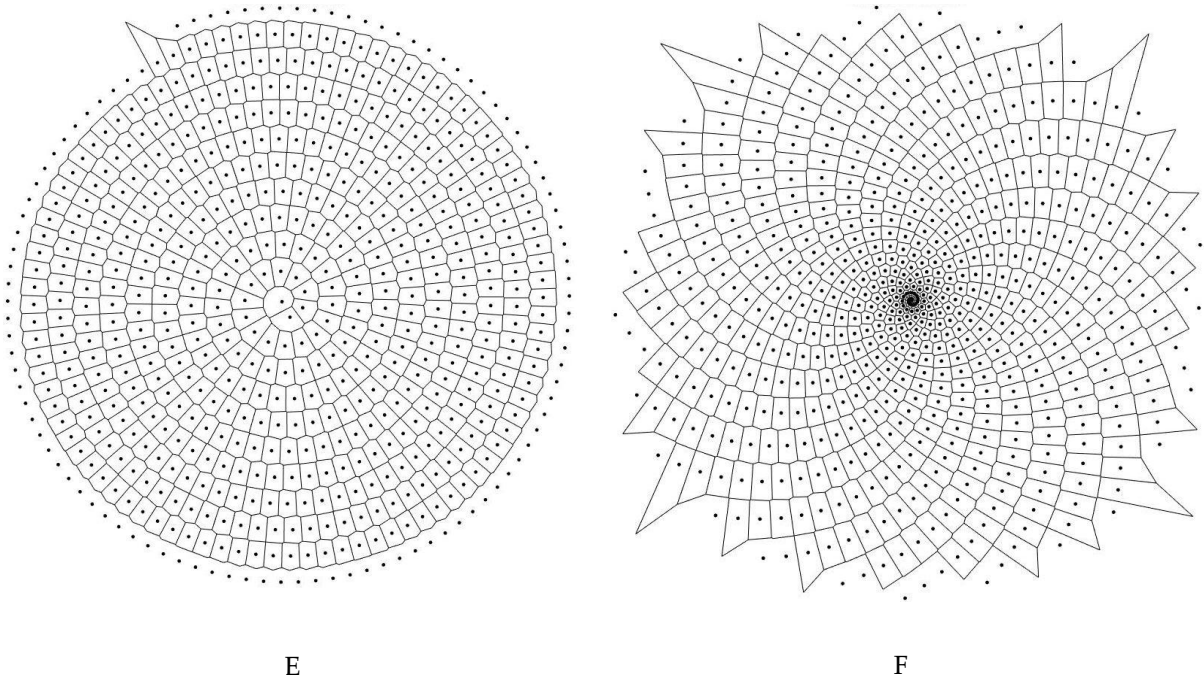

E                                                      F

**Figure 2**. Voronoi diagrams arising from the points located on the AS are shown. **Figures A, C,** and **E** depict Voronoi diagrams for 60, 200 and 600 points, respectively, placed equidistantly on the AS; **Figures B**, **D**, and **F** depict Voronoi diagrams arising from 60 points ($c=0.5$ $d=30$), 200 points ($c=0.5$ $d=100$), and 600 points ($c=0.5$ $d=300$) located on the AS with the linearly increasing distances between them.

Consider first the tessellations where the distance between the seed points and the distance between the turns of the spiral are of the same order of magnitude (in other words the condition $\xi \cong 1$ takes place). Voronoi tessellations arising from the AS with constant and linearly increasing $p$, and the different total points number of 60, 200, and 600 are displayed in **Figure 2**. It is clearly seen, that for the equidistant distribution of points on a spiral (depicted in **Figures 2A, C, E**), the type of pattern does not change with an increase in the total number of points $N$. The configuration of external (boundary) polygons changes with the growth in the total number of points $N$ on the spiral for the patterns with a linear increase of NP distance. It is recognized from **Figures 2B, D, F**, that with an increase in the number of points constituting the spiral $N$, the resulting pattern more and more resembles a sunflower or ordered rows of grains in a corncob, inherent for phyllotaxis.

Now we address the Voronoi entropy calculated for the tessellations presented in **Figure 2**. The distinct prevalence of hexagons is obvious for the patterns with an equidistant distribution of points. This is an immediate consequence of Euler's equation in two dimensions [31]. Let us introduce the number (abbreviated NR) and the area ratios (abbreviated AR) of polygons on the pattern as follows: $NR = \frac{N_e}{N} \times 100\%$; $AR = \frac{A_e}{A} \times 100\%$, where $N_e$ and $A_e$ are the number and area of polygons with $e$ edges respectively and $N$ and $A$ are the total number and area covered by polygons correspondingly.

The area ratio AR of hexagons increases with the growth in the number of points from 73% for 600 points pattern to 94% for 12000 points pattern (as illustrated with **Figures 3 A, B** and **Figures 4 A, B**). For the tessellations based on spirals with linearly increasing distance between NP hexagons also occupied most of the area of the pattern. For example, for the 600 points pattern ($c=0.5$, $d=300$) pattern (**Table 1**) hexagons cover up to 60% of the area, as shown in **Table 1**). The AR of hexagons increases with an increase in the total number of polygons forming the mosaic (**Figure 3 C, D** and **Figure 4 C, D**).

**Table 1**. Polygons distribution characteristics.

| Number of a polygon sides, $e$ | NR, % | | | | AR, % | | | |
|---|---|---|---|---|---|---|---|---|
| | Figure 3A, 600 points ($p=1$, $q=1$) | Figure 3B, 12000 points ($p=1$, $q=1$) | Figure 3C, 600 points ($c=0.5$, $d=300$) | Figure 3D, 12000 points ($b=0$, $c=1$, $d=12000$) | Figure 3A, 600 points ($p=1$, $q=1$) | Figure 3B, 12000 points ($p=1$, $q=1$) | Figure 3C, 600points ($c=0.5$, $d=300$) | Figure 3D, 12000 points ($c=1$, $d=12000$) |
| 3 | 0 | 0 | 0 | 0 | 0 | 0 | 0 | 0 |
| 4 | 0 | 0.01 | 2.15 | 0.01 | 0 | $3.25 \times 10^{-2}$ | 6.522 | $2.63 \times 10^{-6}$ |
| 5 | 13.73 | 3.21 | 16.52 | 4.55 | 13.73 | 3.21 | 18.1123 | 8.5875 |
| 6 | 73.31 | 93.61 | 67.68 | 92.36 | 73.07 | 93.56 | 61.5241 | 89.8156 |
| 7 | 12.77 | 3.17 | 13.46 | 3.09 | 12.85 | 3.18 | 13.8391 | 1.5969 |
| 8 | 0 | 0 | 0 | 0 | 0 | 0 | 0 | 0 |
| 9 | 0 | 0 | 0 | 0 | 0 | 0 | 0 | 0 |
| 10 | 0.19 | 0.01 | 0.18 | 0 | 0.35 | $1.55 \times 10^{-2}$ | $2.41 \times 10^{-3}$ | 0 |

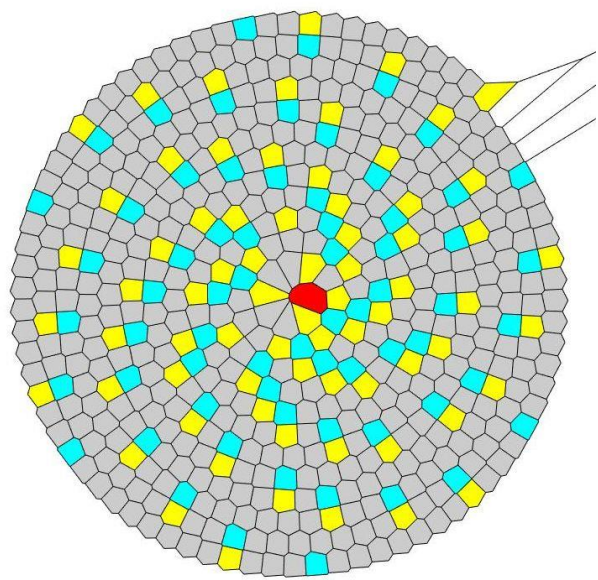

A

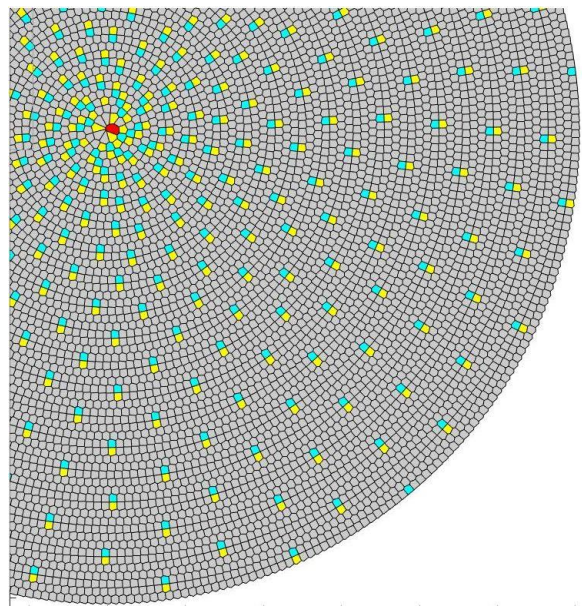

B

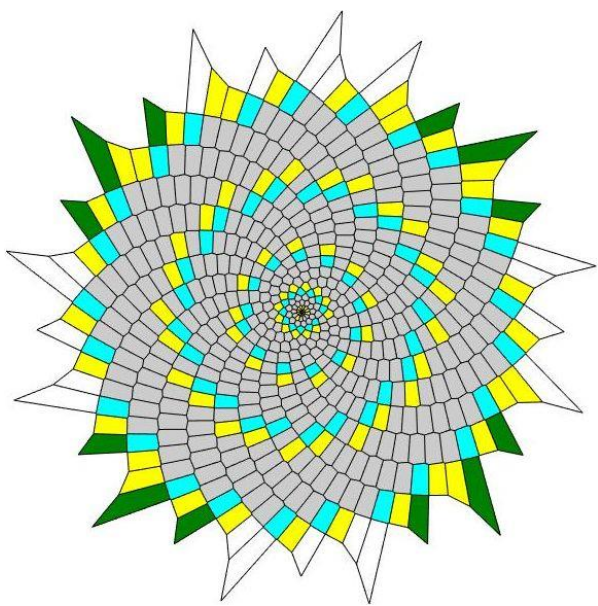

C

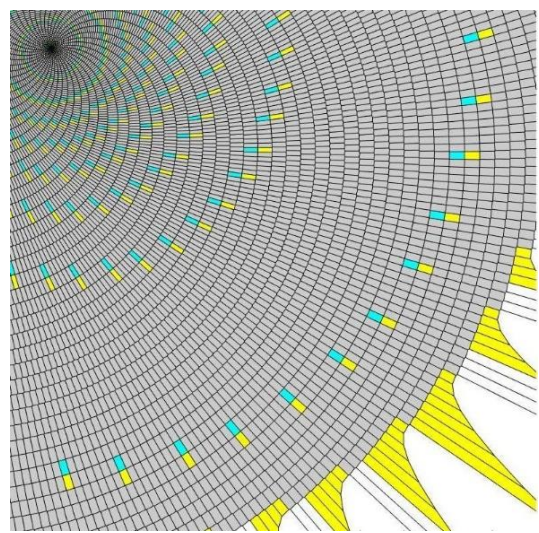

D

**Figure 3**. Voronoi diagrams with colored polygons are shown (grey polygons are hexagons, yellow- pentagons, blue-heptagons, green- tetragons). **A**. A 600-points AS, equidistant points distribution along the curve, (*p*=1, *q*=1). **B**. a 12000-points AS, equidistant points

distribution along the curve, ($p=1$, $q=1$). **C**. A 600-points AS, linear NP distance increase, ($c=0.5$, $d=300$). **D**. A fragment of 12000-points AS, linear NP distance increase, ($c=1$, $d=12000$).

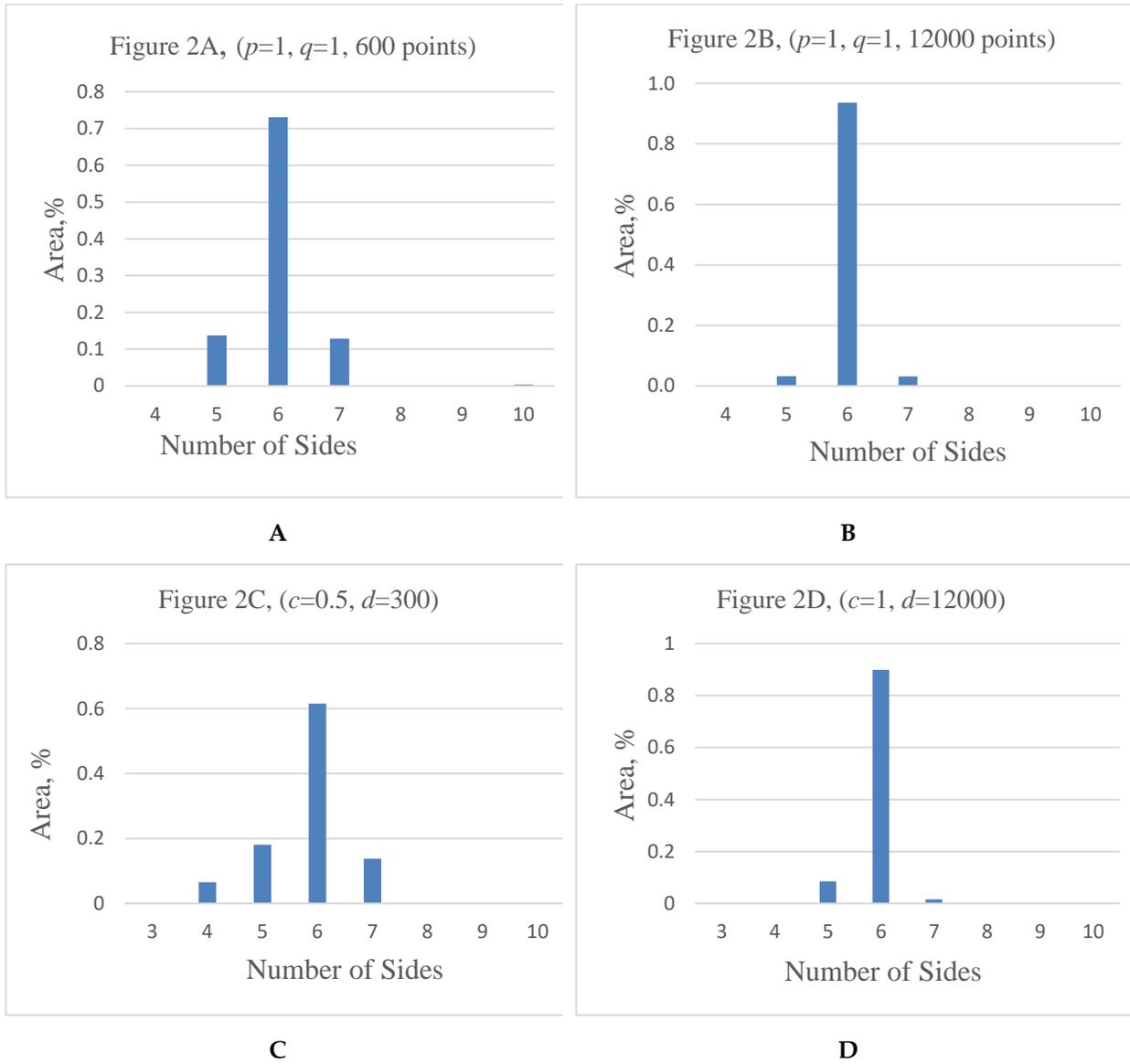

**Figure 4**. The relative area occupied by different types of polygons in the tessellations is shown. **A** and **B** correspond to equidistant location of neighboring points on the spiral. **C** and **D** correspond to a linear increase of the distance between neighboring points on the spiral.

Eq. 1 enabled calculation of the Voronoi entropy for a given pattern. The Voronoi entropy values were obtained for spirals with different density of points and their entire quantity $N$. The

Voronoi entropy $S_{vor}$ depended markedly on the entire quantity of points for the both equidistant and linearly increasing distance points distribution. Consider first patterns arising from the spirals with equidistant NP shown in **Figures 3A, 2B**. It is recognized from **Figure 3B** that actually the Voronoi entropy $S_{vor}$ is built from the contributions supplied by blue heptagons bordering yellow pentagons (which may be called "the defects"). Indeed, the contribution of closely packed hexagons to $S_{vor}$ is negligible. Hence, the value of $S_{vor}$ mainly results from the secondary spiral-like pattern created by the pairs of heptagons bordering pentagons, as illustrated with **Figures 3A,B**.

Similar Voronoi mosaics, arising from the analysis of the Benard-Marangoni cells were reported and discussed in Ref. 27 by Rivier et al. It was noted in ref. 27 that penta- and heptagonal cells represent positive or negative disclinations (corresponding to rotational dislocations, well-known in crystallography) and that they are topologically defined objects which are structurally stable; in other words, they keep their identity under small deformation. Rivier et al. related their appearance to the finite nature of the studied pattern. The defects are a necessary ingredient of the finite mosaic [27]. The reported finite mosaics are necessarily restricted by the origin (the area adjacent to origin is "defected", as shown in **Figures 2A,B**) and boundary points. Thus, the boundary conditions are crucial for the formation of the resulting pattern [27]. In our case the boundary conditions are pre-scribed by the location of the seeds on the AS, given by Eqs. 2-3. Indeed, it is clearly recognized from **Figure 3B** that the larger is the pattern, the smaller are the number and area ratios of the defects. It is also seen from **Figures 3A-D** that pentagons attract heptagons, as reported in Ref. 27.

Now consider the dependence $S_{vor}(N)$. In the case of equidistant location of the seed points along the spiral ($p$=const) the Voronoi entropy is decreasing monotonously with the growth of the total number of points $N$, with an exception of the initial part of the curve $S_{vor}(N)$. For the pattern with $p$=3, $q$=3, we calculated $S_{vor}(10^4) \cong 0.3$, $S_{vor}(2.5 \times 10^4) \cong 0.2$. The initial jump in the curve $S_{vor}(N)$ is due to large density of "defects" (presented as yellow pentagons and blue heptagons in **Figure 3A**) appearing at the initial stage of spiral formation. The density of these "defects" decreases with the growth of the total number of points $N$, leading to the monotonic decrease in the resulting Voronoi entropy of the pattern.

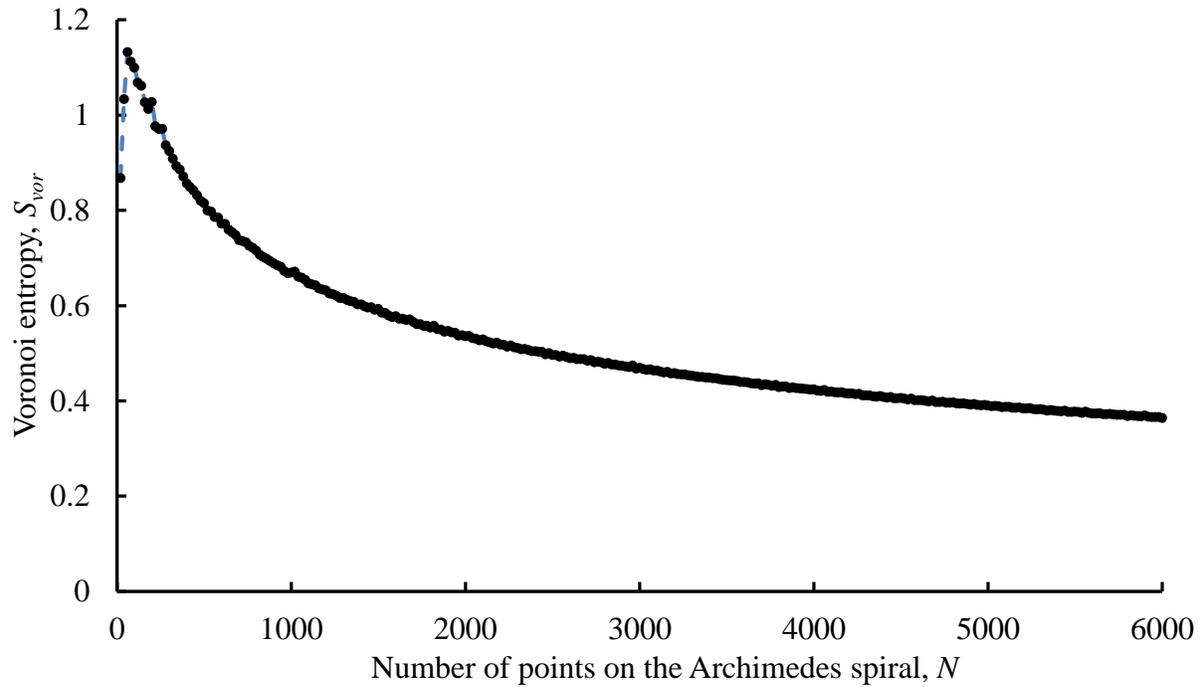

**Figure 5**. The dependence of the Voronoi entropy calculated for patterns arising from equidistant distribution of points located on the AS on the total number of points $N$ is shown ($p=3$, $q=3$, $N$ is changed from 20 to 6000 with a step of 20).

The Voronoi entropy for the patterns based on a spiral with the linear increase of distance between neighboring points demonstrates a saw-like character, depicted in **Figure 6A** while tending to decrease with the growth of the entire number of points $N$. The saw-like behavior of the curve $S_{vor}(N)$ is reasonably explained as follows: the bordering heptagons and pentagons form the ring-like secondary pattern, contributing markedly to the $S_{vor}(N)$. Appearance of these rings (introducing geometrical disorder into the pattern) increases the value of $S_{vor}(N)$, resulting in the saw-like dependence of the function $S_{vor}(N)$.

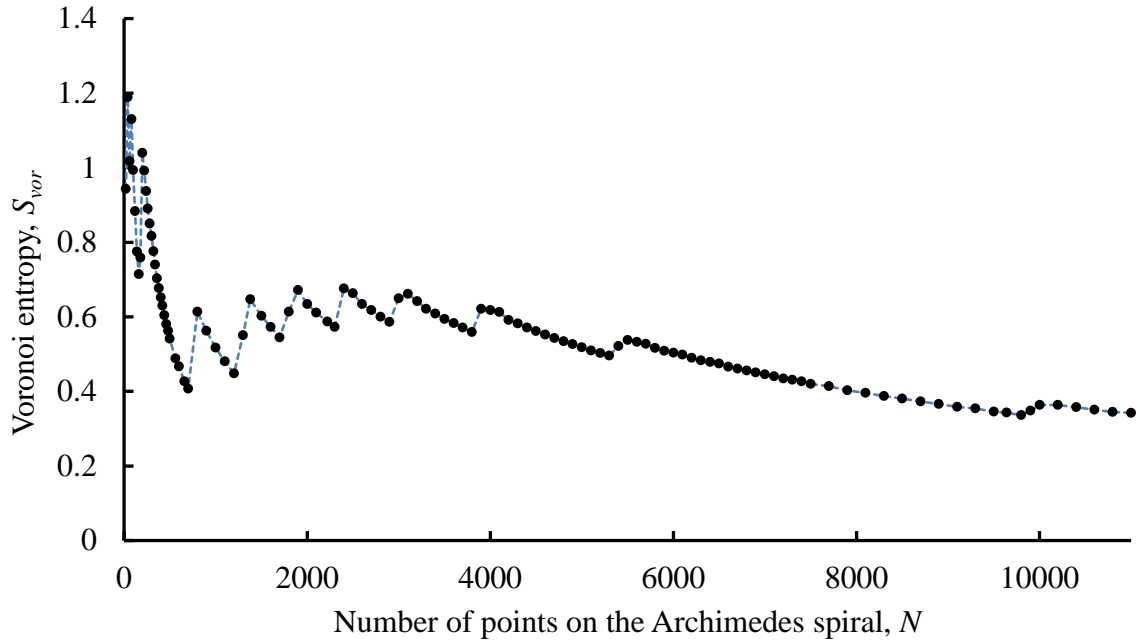

A

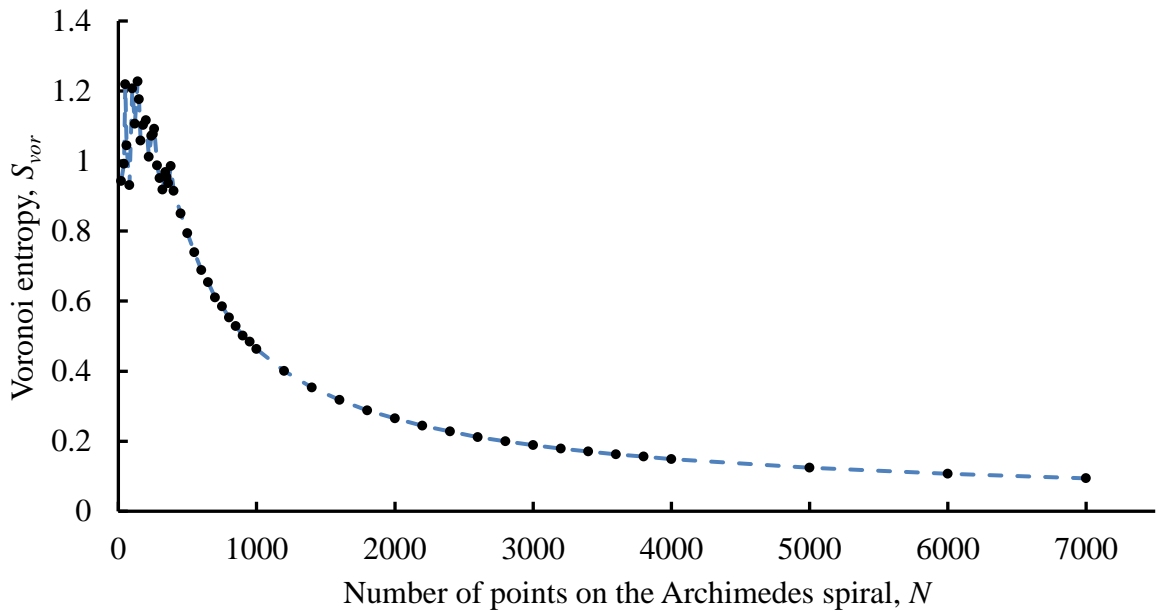

B

**Figure 6**. Voronoi entropy change with increasing number of points $S_{vor}(N)$ for the AS with the linearly increasing distance between neighboring points is shown. **A**. Parameters of the pattern: $c=1$, $d$ varies from 20 to 11000 with a step of 20. **B**. Parameters of the pattern: $c=20$, $d$ varies from 20 to 400 with a step of 20, from 400 to 1000 with a step of 50 and from 1000 to 7000 with a step of 200.

We relate the origin of these picks' appearance to the appearance of irregularities on Voronoi diagrams. Two kinds of irregularities (defects) inherent for AS-inspired Voronoi diagrams should be distinguished, the first of which is the fringe of a Voronoi pattern formed by open (incomplete) polygons. The fringes-effect on the Voronoi entropy is essential for patterns consisting of a small number of points (polygons, respectively). The second type of irregularities is represented by the aforementioned "defected areas" of a pattern (colored with blue and yellow in **Figure 3D**. These irregularities appear as circles of pentagons bordering heptagons recognized on the background filled with hexagons [27]. The well- ordered nature of the defected areas is noteworthy. With an increase in the number of points, the distance between irregularity circles is growing, as shown in **Figure 3B**. The growth of the distance between blue/yellow circles leads to a consequent decrease in the Voronoi entropy of the entire pattern.

It is also noteworthy that in the case of linearly increasing distance between the adjacent seed points it is possible to select values of parameters of $c$ and $d$ (for example consider the case of $c=20$, $d=40000$) resulting in the disappearance of irregularities built of pentagons and heptagons (this occurs for $N>N^*$, where $N^*$ is the threshold value of points, corresponding to two or three central rings). In this case, the value of Voronoi entropy falls faster and asymptotically tends to zero. When $N\to\infty$ the role of the central area of the spiral becomes negligible and $S_{vor} \to (0+ \frac{1}{\infty})$. Consequently, the Voronoi entropy $S_{vor}(N)$ does not show the saw-like behavior, when the "defected circles" are absent, as shown in **Figure 6B.** This possibility to fill a plane with cells of equal size is of a primary importance for phyllotaxis (leaf or floret arrangement) and decorative arts [27].

### 3. The Aboav and Lewis laws for the patterns inspired by the Archimedes Spiral

The Aboav law validity for the mosaics generated by AS was verified. The Aboav law states that the mean number of sides of polygons (labeled $m_n$) bordering the polygon with $n$-edges is given by [1,5,24,30]:

$$m_n = 5 + \frac{8}{n} \qquad (4)$$

In other words, the few-edged cells have a remarkable tendency to be in contact with many-edged cells and vice versa. The critics, derivation and consequences of the Aboav law are discussed in Refs [1,24,30,31]. The values of $m_n$ were calculated for the pattern with linear increase of distance between NP ($c=1$, $d=500$) and for the pattern with equidistant points distribution ($p=3$,

$q=3$, $N=500$). All of the obtained $m_n$ values were lower than those calculated with the formula (4). Thus, the Aboav law, supplied by Eq. 4 does not work for the Voronoi patterns inspired by AS. We relate this observation for a non-random distribution of points on the studied mosaics.

On the other hand, the mean values of $m_n$ for the AS-inspired patterns tend to decrease with an increase in the number of sides of the corresponding polygon, as it qualitatively predicted by the Aboav law [1,5,24,30].

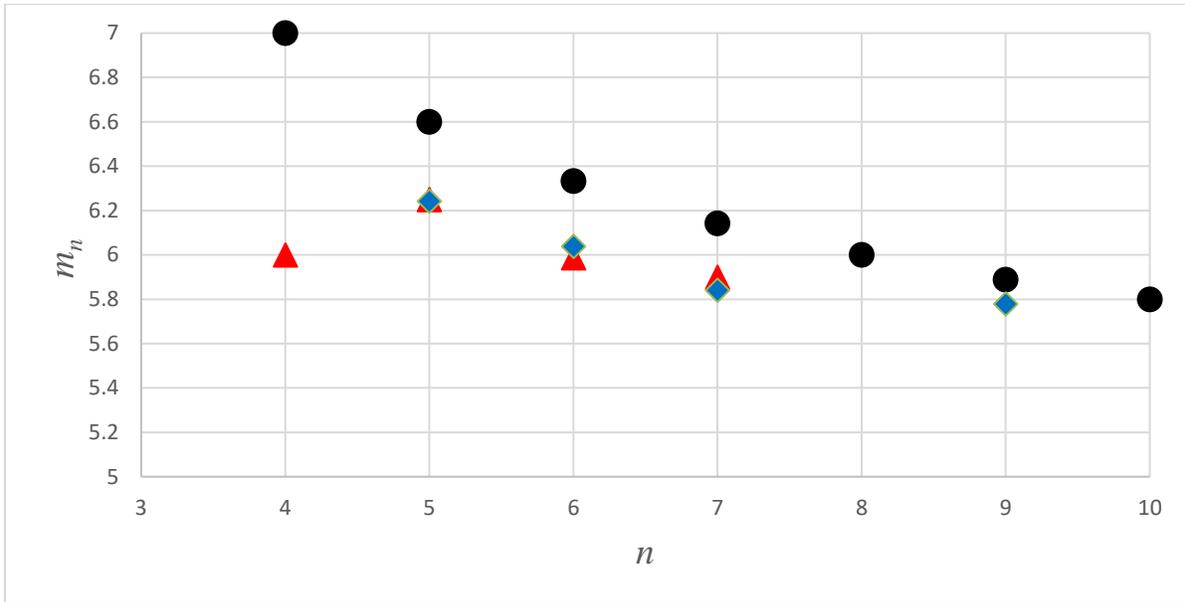

**Figure 7.** Variation of $m_n$ with the number of polygon sides $n$ is shown. Black circles – values of $m_n$ calculated with the Aboav law (eq. (4)); red triangles are the average values of $m_n$, obtained the pattern with linearly increasing NP distance ($c=1$, $d=500$), blue diamonds are average values of $m_n$, obtained for the pattern with equidistant points distribution ($p=3$, $q=3$, $N=500$).

Another important statistical law, established for 2D patterns is the Lewis law, reported for natural and artificial patterns [5,17,18]. The Lewis law predicts a linear relationship between the average area of a typical $n$-cell, $A_n$ and $n$ in a random pattern:

$$A_n = \alpha(n - 2), \qquad (5)$$

where $\alpha$ is a proportionality constant, which meaning and precise value can be found in Ref. [30]. The Lewis law quite expectably does not work for Mosaics generated by Archimedes Spirals. For the patterns based on equidistant NP distribution, polygon areas have a constant mean value of

9.0±0.01 mm² all over pattern. In the case of patterns with linearly increasing NP distance, the areas of polygons on spiral coils are growing with a distance from the origin of a spiral.

4. **Patterns generated by Archimedes Spirals and the maximal Voronoi Entropy.**

Consider patterns with an equidistant distribution of points, in which the NP distance $p$ is much greater (an order of magnitude) than the distance between turns of a spiral $q$. When $p >> q$ and correspondingly $\xi \gg 1$ takes place, we assume $\varphi = k\pi + \Delta\varphi$, where $k$ is a positive integer. Two examples of such patterns are shown in **Figures 8A,B**. Such patterns contain more kinds of polygons than patterns, where $p \cong q$ and correspondingly $\xi \cong 1$ is adopted. Eight types of polygons constituting these mosaics were registered, when $\xi \gg 1$. It should be emphasized that hexagons do not prevail when $\xi \gg 1$ takes place (see **Figures 7 A,B,E,F**).

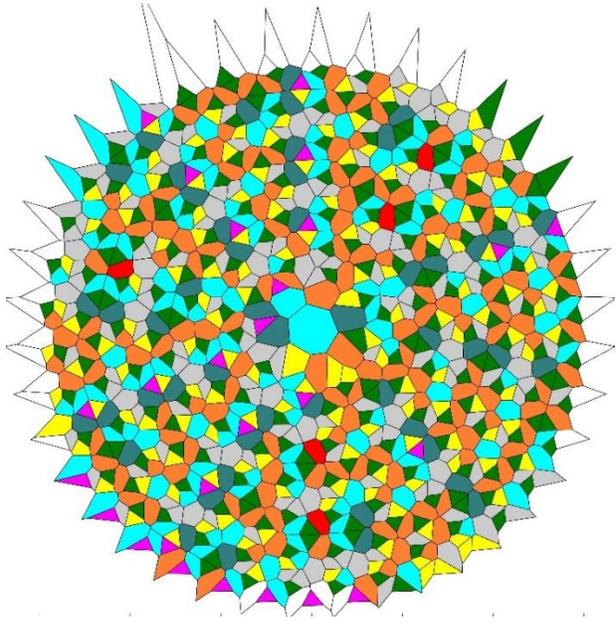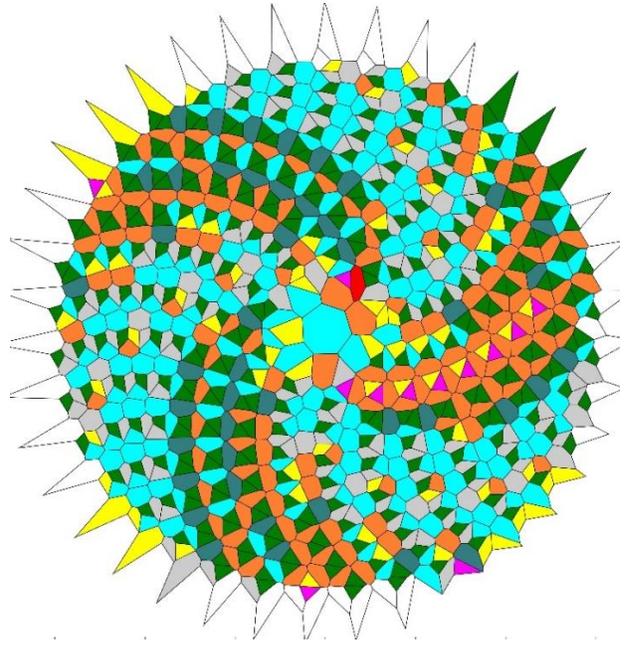

| A | B |

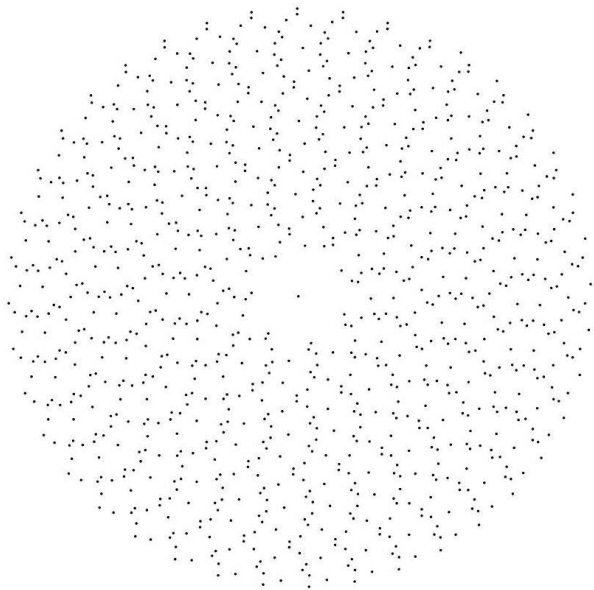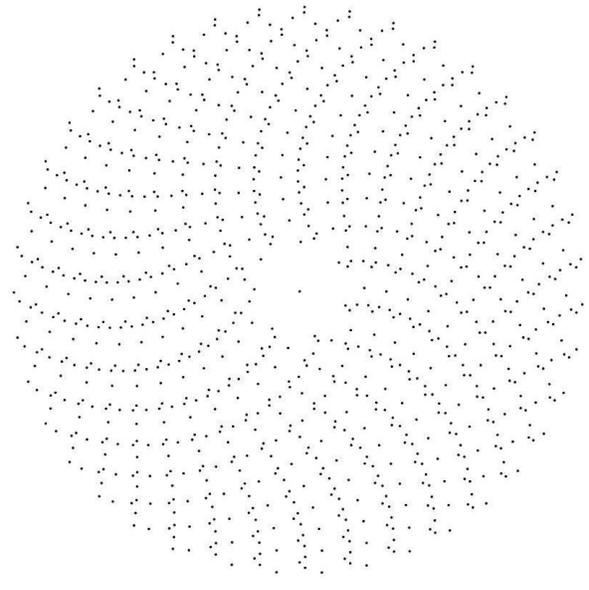

| C | D |

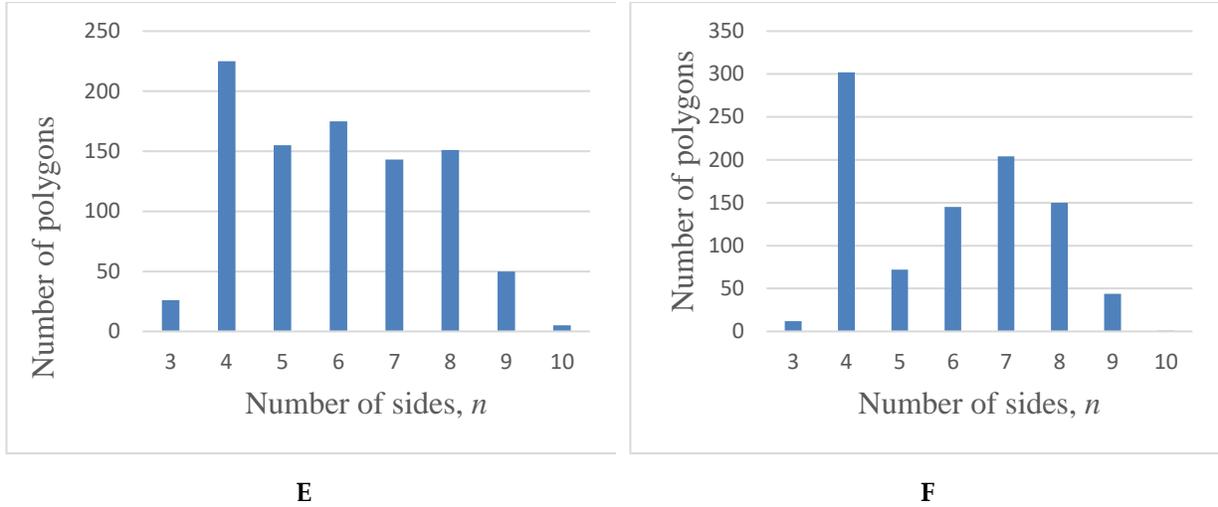

**Figure 8.** Patterns containing $N=400$ points for which $p \gg q$ takes place are shown. **A**. $p = 24.52$, $q=3$, $N=1000$, $S_{vor} = 1.825$; **B**. $p = 25$, $q=3$, $N=1000$, $S_{vor}=1.688$. Color mapping: magenta polygons are triangles, green are tetragons, yellow are pentagons, grey are hexagons, blue are heptagons, brown are octagons, deep-green are nonagons, red are decagons. Figures **C** and **D** depict the location of the seed points. Figures **E** and **F** depict the distribution of polygon kinds in the patterns shown in **Figures 8A** and **8B** correspondingly.

The patterns, where $p \gg q$, are interesting because while not being random, they show high Voronoi entropy values which are close to the value of 1.71, which is considered as maximum inherent for a random pattern [20]. **Figure 8A** depicts the pattern demonstrating the Voronoi Entropy $S_{vor} = 1.825$. The maximal value of the Voronoi entropy $S_{vor} = 1.888$ was registered for the pattern ($p=24.6131$, $q=3$, $N=80$), presented in **Figure 9A**, which is markedly higher than the value reported for random patterns [20,21]. The value of the Voronoi entropy may be even extended to larger values. Consider the pattern arising from the seven-fold $X$ and sevenfold $Y$ translation of the pattern shown in **Figure 9A**. Such a procedure gave rise to the Voronoi diagram shown in **Figure 9B** built from 8 types of polygons and characterized by the Voronoi entropy $S_{vor} = 2.08$, which is much larger than that, established for random point patterns [20, 21].

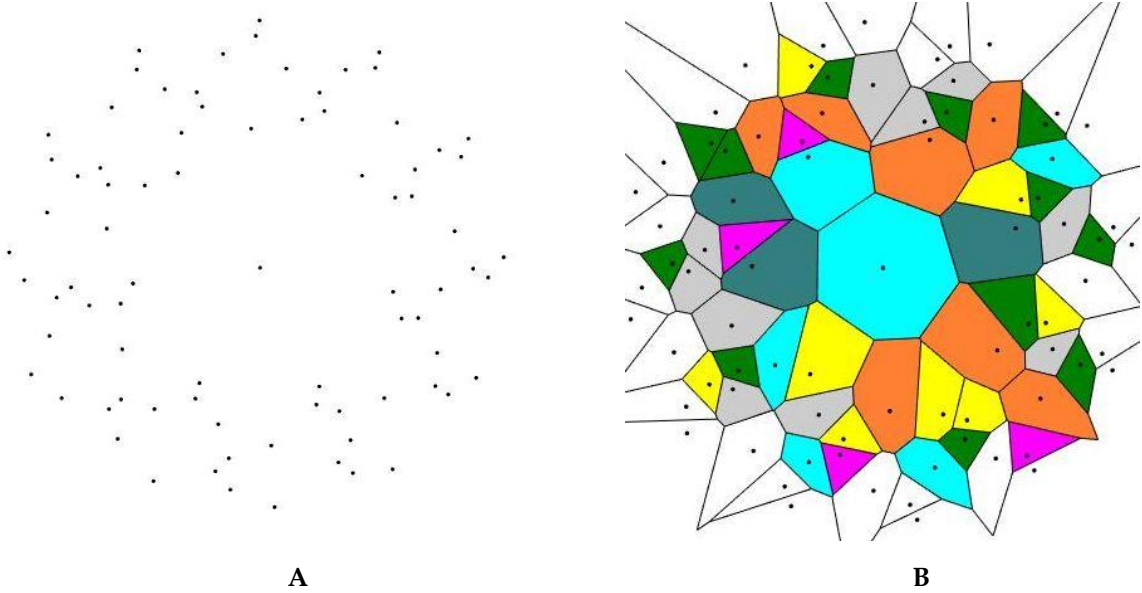

**Figure 9**. An 80 points pattern (*p*=24.6131, *q*=3, *N*=80) (**A**) giving rise to the Voronoi tessellation (**B**) including 7 types of polygons demonstrating the Voronoi Entropy $S_{vor}$=1.8878 is shown. Color mapping: magenta polygons are triangles, green are tetragons, yellow are pentagons, grey are hexagons, blue are heptagons, brown are octagons, deep-green are nonagons.

This finding poses the following fundamental question: it is well-accepted that the Voronoi Entropy quantifies ordering in 2D patterns [2,10,19,29,31]. It is reasonable to conjecture that the maximal disorder corresponds to the random distribution of seeds points over the plane. Hence, the maximal possible Voronoi Entropy is expected for the random distribution of points. At the same time patterns depicted in **Figure 9** are definitely ordered, however possessing the Voronoi Entropy markedly higher that $S_{vor} = 1.71$, established for the random patterns [20,21]. How is this possible? This question calls for additional theoretical insights. Actually, it is well-known that the Voronoi entropy may be larger than 1.71. The maximal value of the Voronoi entropy for the mosaics built from *n* kinds polygons corresponds for a pattern at which equipartition of polygons takes place (their appear $\frac{1}{n}$ of all kinds of n-polygons in the pattern) [28]. In this case the maximal value of the Voronoi entropy is given by $S_{vor}^{max} = ln(n)$ [28], and it obviously may be larger than $S_{vor} = 1.71$, inherent for a random 2D pattern [20,21].

We also have checked the validity of the Lewis law [5,17,18] for the patterns characterized by the $p \gg q$ interrelation (see **Figure 10**). In this case the dependence $An(n)$ demonstrates the linear part, predicted by the Lewis law (see eq. 5), as shown in **Figure 10**.

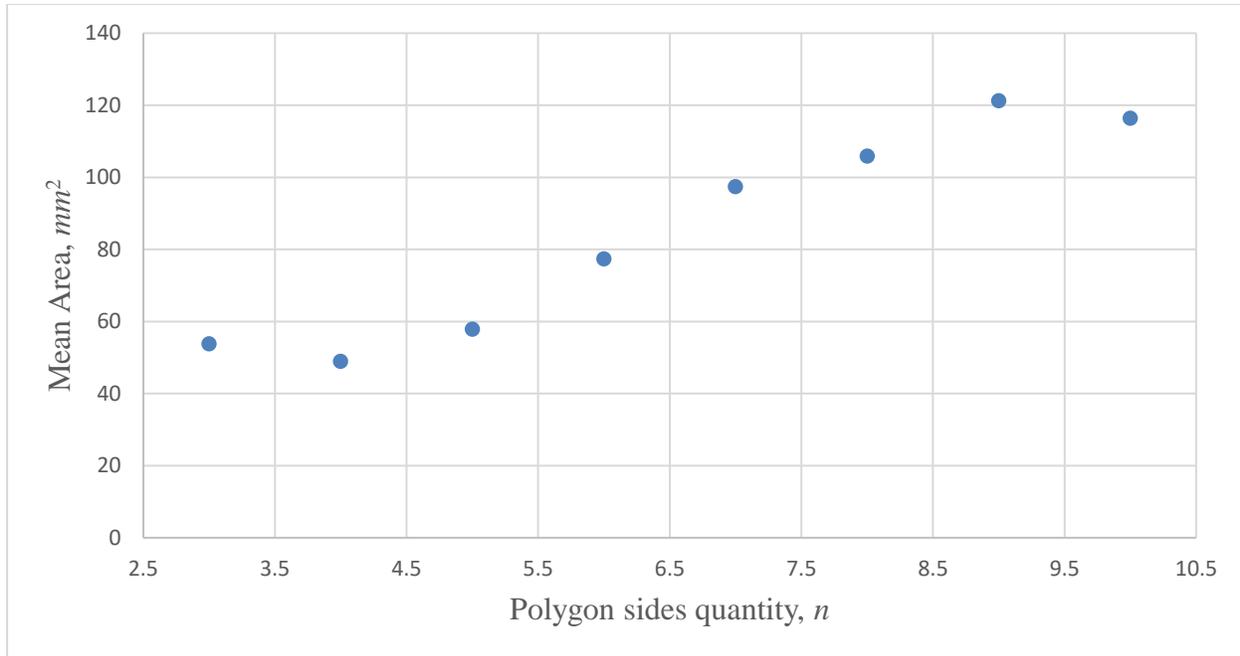

**Figure 10**. Variation of $A_n$ on number of polygon sides $n$ for the pattern depicted in **Figure 8A** is shown.

5. **Voronoi diagrams arising from Archimedes Spirals demonstrating an aesthetic appeal**.

Voronoi diagrams depicted in **Figures 11A-D** demonstrate definite aesthetic appeal. The aesthetic attractiveness of the AS was already known to Neolithic artists [35]. Remarkably AS holds its aesthetic appeal in the XXI century [3]. The very question is: why the spirals demonstrate obvious aesthetic appeal? Obviously, the true answer could not be covered by physics and mathematics only, but also enrooted in psychology [15]. We allow ourselves to put forward the following hypothesis: the aesthetic appeal and abundance of spirals in nature is not least related to their simplicity and self-similarity (the equation describing AS is one of the simplest possible ones). The simplicity and self-similarity are not synonymous but bordering notions [16,22,26,33]. Mathematicians have customarily regarded a proof as beautiful if it conformed to the classical ideals of brevity and simplicity [33]. Similarly, Michael Atiyah claims that "elegance is more or less synonymous with simplicity" [16]. We are quoting from Ref. 26 "The mathematical concept of similarity holds one of the keys to understanding the processes of growth in the natural world. As a member of a species grows to maturity, it generally transforms in such a way that its parts maintain approximately the same proportion with respect to each other, and this is probably a reason why nature is often constrained to exhibit self-similar spiral growth". AS and AS-inspired Voronoi mosaics exemplify simple, self-similar structures. This at least partially explains their aesthetic appeal. The concept of "beauty as simplicity" was strongly criticized recently [16]; thus, additional insights into understanding of the aesthetic appeal of spiral-inspired patterns are necessary. Of course, beauty of mathematical proofs and images (pictures) may be very differently psychologically rooted; however, it is also possible that the concept of "self-similarity" is common for our aesthetic estimation of a broad range of heterogeneous phenomena.

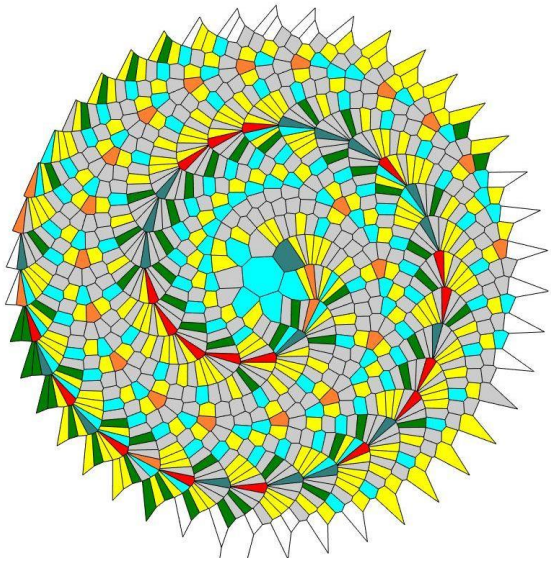 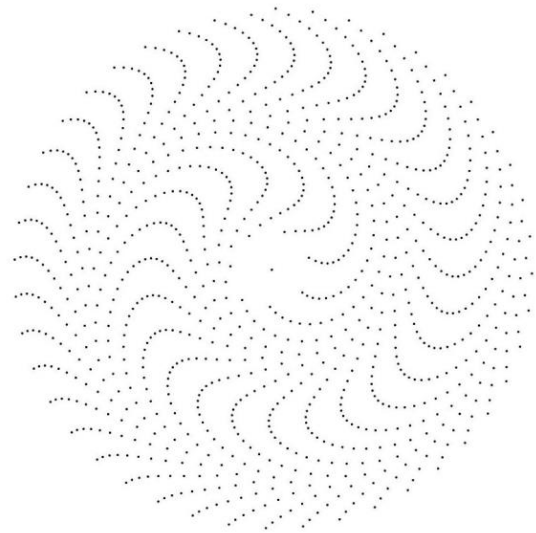

**A**

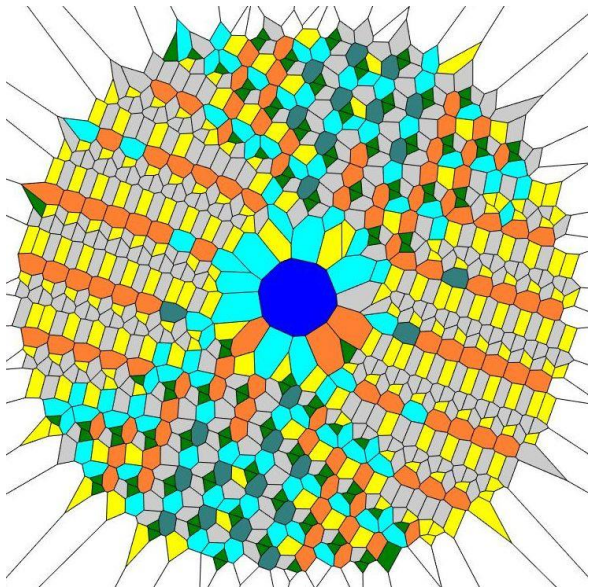 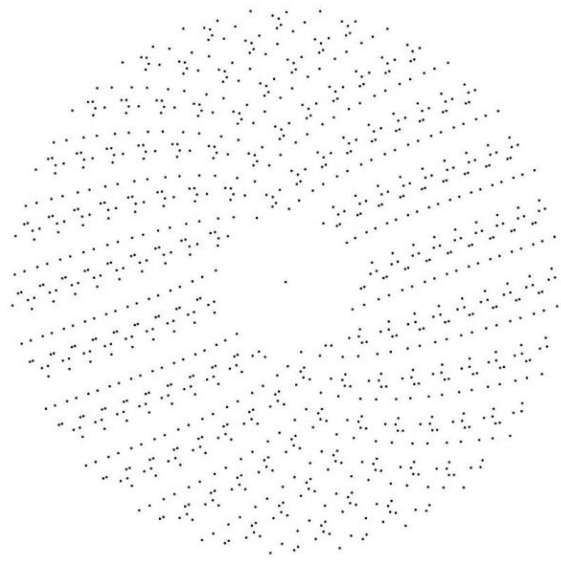

**B**

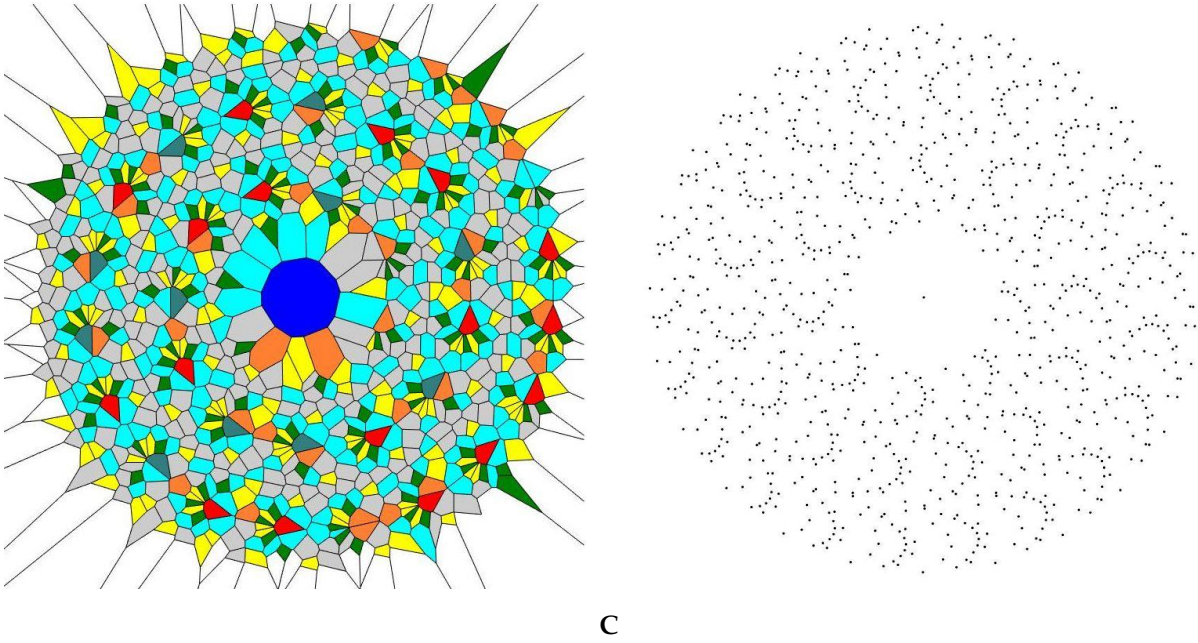

**Figure 11.** Voronoi tessellations demonstrating definite aesthetic appeal are shown. The tessellations were built on equidistant points distribution on AS for which *p>>q* takes place. **A**. $p = 19.75$, $q=3$. **B**. $p = 66$, $q=3$. **C**. $p = 72.5$, $q=3$. Color mapping: green polygons are tetragons, yellow are pentagons, grey are hexagons, blue are heptagons, brown are octagons, teal are nonagons, red are decagons.

### Conclusions

We conclude that the Voronoi diagrams generated by seed points located on the Archimedes Spirals demonstrate non-trivial mathematical properties and aesthetic attraction. Equidistant seed point distribution and points separated by linearly increasing distance generated very different Voronoi diagrams. Voronoi entropy calculated for the equidistant seed points located on the Archimedes Spiral decreased monotonously with the increase in the number of seeds. The Voronoi entropy calculated for points separated by linearly increasing distance demonstrated a saw-like behavior. It is possible to fill a plane with Voronoi mosaics built from a cells of equal size which is of a primary importance for phyllotaxis and decorative arts.

The properties of the Voronoi tessellation to a much extent are governed by the parameter $\xi = \frac{p}{q}$ where *p* and *q* are the distance between the points neighboring along the spiral and the separation between the coils of the spiral, respectively. When the condition $\xi \cong 1$ is assumed, hexagons dominate in the mosaic; whereas, eight types of polygons were registered, when the condition $\xi \gg 1$ was prescribed. For the patterns characterized by $\xi \gg 1$ the ordered patterns were revealed demonstrating the Voronoi entropy markedly larger than that of 1.71, reported for the random distribution of points [20,21].

Archimedes Spirals generate Voronoi diagrams enabling to fill a plane with equal size cells. This possibility is of a primary importance for phyllotaxis (leaf or floret arrangement) [27]. The Aboav and Lewis laws generally do not hold for the Voronoi mosaics generated by the Archimedes Spirals. We explain this observation by the non-random distribution of seed points inherent for studied patterns. The Voronoi mosaics inspired by the Archimedes Spirals demonstrate definite aesthetic appeal [3,15,35]. We relate at least partially the aesthetic attraction of the reported mosaics to their simplicity and self-similarity [22,33]. In our future work we plan to consider symmetry considerations applied to the analysis of Voronoi diagrams inspired by the Archimedes Spirals [27].

## Acknowledgments

Dr. Mark Frenkel acknowledges partial support from the Israel Ministry of Immigrant Absorption. The authors are indebted to the anonymous reviewer for a fruitful reviewing of the manuscript.

## Conflict of Interests

The authors declare no conflict of interests.

## Appendix A

Coordinates of points on AS in rectangular coordinate system can by defined by following equations:

$$x = r \cdot \cos(\varphi), \qquad (2)$$
$$y = r \cdot \sin(\varphi).$$

In the case of a linear increase of distance between NP we increased parameters $r$ and $\varphi$ discretely from $b$ (which was equal to 0 for of all studied patterns) to $d$ with a step of $c$. Thus, the magnitude of the $r$ increment was constant and equal to $c$ in each of subsequent steps (See **Figure A1**). The magnitude of the angle $\varphi$ increment was also constant and equal to $c$ (See **Figure A2**). These parameters were matrixes of same size varied in the same way. Therefore, in our MATLAB routine parameters $r$ and $\varphi$ were substituted by the single parameter denoted $t$ as follows:

$$x = t\{b,c,d\} \cdot \cos(t\ \{b,c,d\}),$$
$$y = t\{b,c,d\} \cdot \sin(t\ \{b,c,d\}), \qquad (3)$$

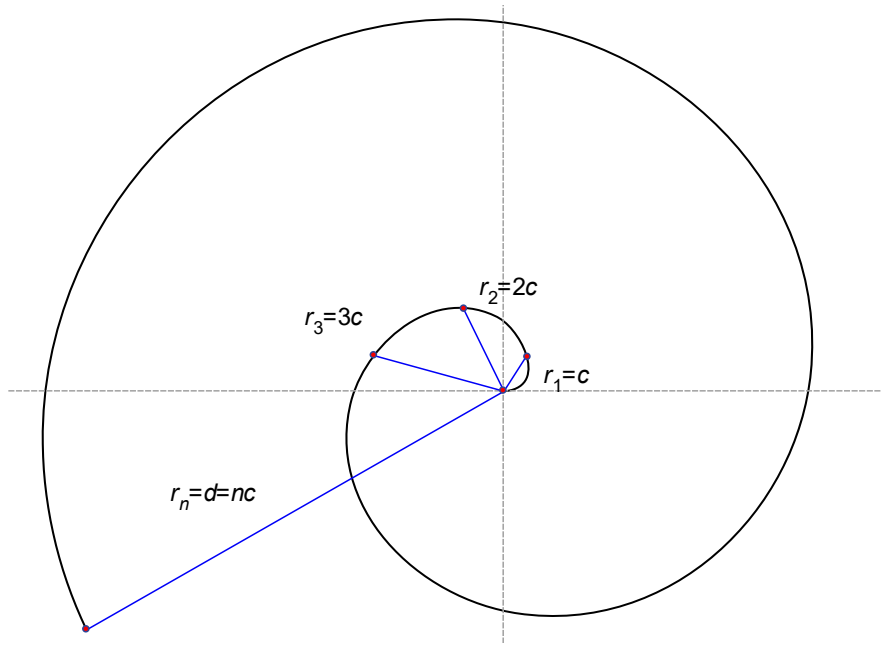

**Figure A1**. The change in $r$ is shown.

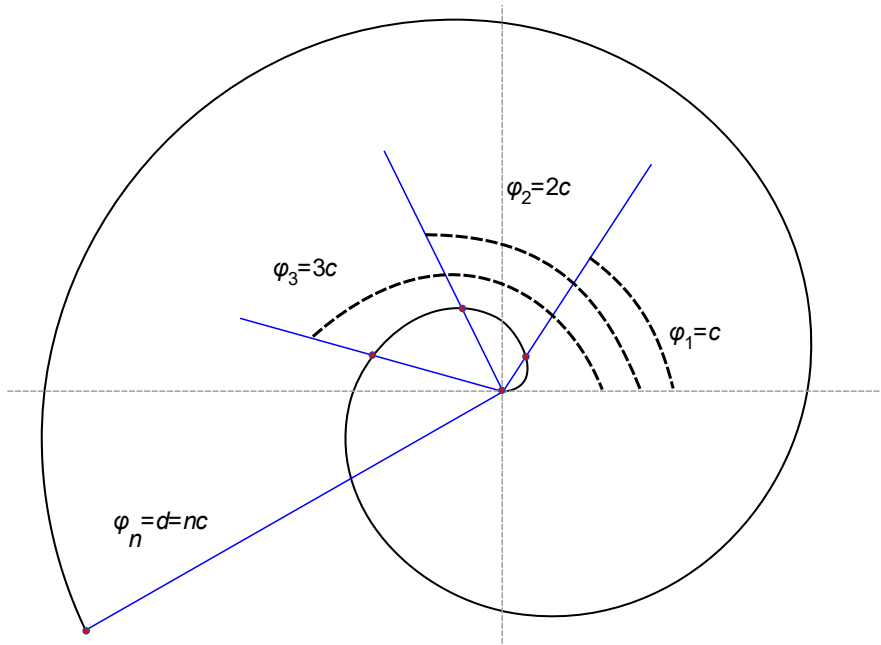

**Figure A2**. The change in $\varphi$ is shown.

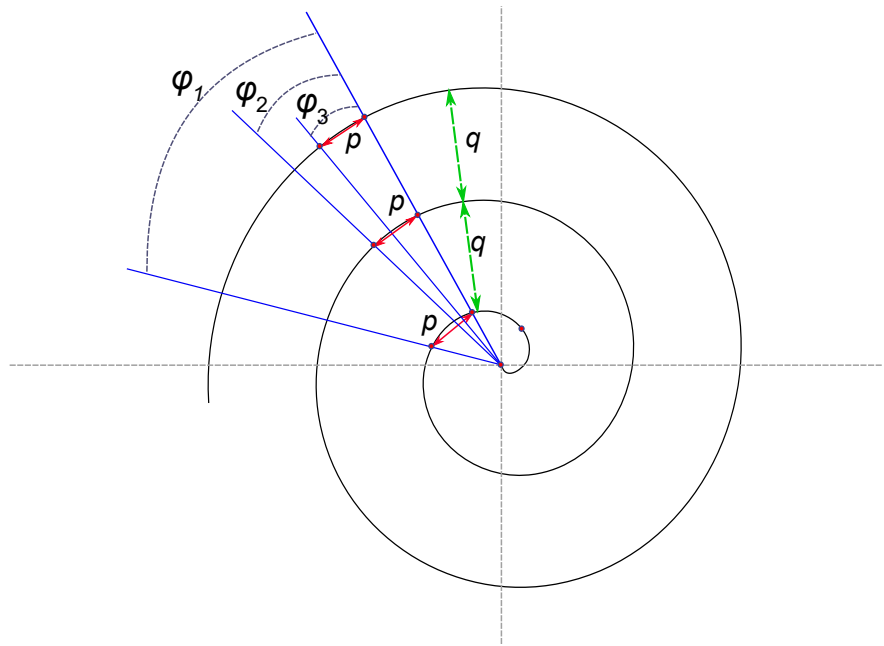

**Figure A3**. Parameters $p$ and $q$ of the points' distribution along the Archimedean spiral are shown.

In the case of equidistant points' distribution prescribed along the Archimedean spiral, we increased $r$ and $\varphi$ in a way to keep parameters $p$ and $q$ constant. To achieve this, we had to reduce the increment of the angle in each step (See **Figure A3**).